\documentclass[12pt,leqno]{article}
\usepackage{latexsym}
\usepackage{amsmath,amssymb}
\usepackage{mathrsfs}
\usepackage{amsthm}
\usepackage{amsfonts}
\numberwithin{equation}{section}
\newtheorem{definition}{Definition}[section]
\newtheorem{theorem}{Theorem}[section]
\newtheorem{lemma}{Lemma}[section]

\newtheorem{proposition}{Proposition}[section]

\newtheorem{remark}{Remark} [section]
\newtheorem{corollary}{Corollary}[section]

\begin{document}
\author{Hanzhong Wu \\School of Mathematical Sciences\\ Fudan University, Shanghai 200433, China\\Email:
hzwu@fudan.edu.cn}
\date{}
\title{A Novel Variational Principle arising from Electromagnetism} \maketitle
\noindent{{\bf Abstract. } Analyzing one example of LC circuit in
\cite{MA}, show its Lagrange problem only have other type critical
points except for minimum type and maximum type under many
circumstances. One novel variational principle is established
instead of Pontryagin maximum principle or other extremal principles
to be suitable for all types of critical points in nonlinear LC
circuits. The generalized Euler-Lagrange equation of new form is
derived. The canonical Hamiltonian systems of description are also
obtained under the Legendre transformation, instead of the
generalized type of Hamiltonian systems. This approach is not only
very simple in theory but also convenient in applications and
applicable for nonlinear LC circuits with arbitrary topology and
other additional integral constraints.
}\\

\section{Introduction}

The Lagrangian and Hamiltonian formulation of nonlinear
inductor-capacitor circuits (LC circuits) has been considered by
[3--8], \cite{SJK} and the many references incited within. van der
Schaft, Maschke and coworkers \cite{MVB} and Bloch and Crouch
\cite{BC} etc. established the Hamiltonian modeling by utilizing the
constant Dirac structure of circuits.

Many authors considered the variational approaches to nonlinear LC
circuits. Based on the dual extremum principle in \cite{NS}, Kwatny,
Massimo and Bahar \cite{KMB} realized the Lagrangian modeling.
Recently, Moreau and Aeyels \cite{MA} obtained the generalized
Euler-Lagrange equations and the generalized Hamiltonian system
description after applying Pontryagin maximum principle. At the same
time, Scherpen, Jeltsema and Klaassens \cite{SJK} established the
Lagrangian modeling of nonlinear LC circuits based on the
constrained variational principle from holonomic mechanics, see e.g.
\cite{Va, Bl} etc. Then Clemente-Gallardo and Scherpen \cite{CS}
considered the relation between Lagrangian and Hamiltonian
formalisms of nonlinear LC circuits via Lie algebroid.

The Lagrangian functional in LC circuits has its own rich
properties. As pointed out in \cite{MA} (see also \cite{KMB}), the
generalized velocities are not simply the derivatives of the
generalized coordinates. In other words as in \cite {CS}, it lacks
of kinetic terms for the capacitors and the potential terms for
inductors. Utlizing the notion of tree and cotree in \cite{SR},
\cite{MA} developed one method to consider nonlinear LC circuits
with any topology (including with excess elements or without excess
elements) through the Lagrangian functional
\begin{eqnarray}
J(u)=\int_{t_0}^{t_1}L(x,u)\,dt,
\end{eqnarray}
subject to the dynamics (from Kirchhoff's current law)
\begin{eqnarray}
\dot{x}=Au,
\end{eqnarray}
for some matrix $A$.

As shown in section 2, this Lagrange problem only have other type
critical points except for minimum type and maximum type under many
circumstances.

In this paper, we will establish a variational principle for this
Lagrange problems. This principle can derive the generalized
Euler-Lagrange equation of new type to describe critical points of
all types. Meanwhile, the Hamiltonian and the canonical Hamiltonian
systems formulation will be given uniformly to describe the critical
points of all types under the generalized Legendre transformation.

One of the advantages of this generalized Euler-Lagrange equation
formulation is that we can very easily derive the canonical
Hamiltonian systems under the generalized Legendre transformation.
In this way, the energy function can be explicitly constructed,
which is very important especially in applications. It should be
pointed out the notion of the generalized Legendre transformation is
adapted from the ideas of H. J. Sussmann and J. C. Willems
\cite{SW}.

The second advantage of this generalized Euler-Lagrange equation
formulation is that it clearly indicate the distinction between the
problems without additional constraints and those with constraints,
especially the terminal state constraints. Meanwhile, it will be
clearly shown in this canonical Hamiltonian system formulations how
the constraints influence the constructions of the Hamiltonian
functions -- the energy functions.

The third advantage is that, in the applications to nonlinear LC
circuits, we will not encounter the technical difficulty to consider
the abnormal cases which unavoidably arising in applying Pontryagin
maximum principle with additional constraints such as terminal state
constraints. Meanwhile, we will also not encounter the complex
calculations of pseudo-inverses of matrixes and Lagrangian
multipliers, which involve in applying the dual extremal principles
of \cite{NS}.

This paper is organized as follows. In section 2, we analyze the LC
examples considered in \cite{MA}. It will be shown that the Lagrange
problem with constraints has other type critical points except for
both minimum and maximum type critical points under many
circumstances; and then put forward a new type of variational
problem instead of minimizing problems of the Lagrange functional.
In section 3, we will establish one variational principle to derive
the generalized Euler-Lagrange equations of new type to describe the
critical points of all types. Some illustrative examples will be
given. In section 4, we will derive the canonical Hamiltonian
systems to describe the critical points of all types under the
generalized Legendre transformation. Some illustrative examples will
also be given. Last, an appendix will be attached enclosed within
the proofs of the results in section 2.

\section{A New Variation Problem arising from Inductor-Capacitor Circuits}

L. Moreau and D. Aeyels \cite{MA} comprehensively considered the
dynamic equation of one LC circuit as illustrative examples (see
Examples 1, 2, 3 and 4 in \cite{MA}).

The associated Lagrange functional consists of magnetic coenergy (of
inductors) minus electric energy (of capacitors)
\begin{eqnarray}\label{ma}
\begin{array}{ll}
J(i_3,i_5,i_6)=&\int_{t_0}^{t_1}[\frac{1}{2}L_3i_3^2+\frac{1}{2}L_4(i_3-i_5-i_6)^2+\frac{1}{2}L_5i_5^2+\frac{1}{2}L_6i_6^2]\,dt\\
&-\int_{t_0}^{t_1}[\frac{1}{2C_1}q_1^2+\frac{1}{2C_2}q_2^2]\,dt,
\end{array}
\end{eqnarray}
where $q_1$ and $q_2$ are described by the dynamics (from
Kirchhoff's current law)
\begin{eqnarray}\label{ma1}
\dot{q_1}=i_3,\qquad \dot{q_2}=i_5+i_6.
\end{eqnarray}
In addition, the following other integral constraints were imposed
in \cite{MA}:
\begin{eqnarray}\label{ma2}
\int_{t_0}^{t_1}i_3dt=\lambda_3,\quad\int_{t_0}^{t_1}i_5dt=\lambda_5,\quad\int_{t_0}^{t_1}i_6dt=\lambda_6.
\end{eqnarray}

In \eqref{ma}-\eqref{ma2}, $C_1>0$ and $C_2>0$ are capacitance,
$L_3>0$, $L_4>0$, $L_5>0$ and $L_6>0$ are inductance, $i_3$, $i_5$
and $i_6$ are the currents. For more information in detailed, please
see Examples 1, 2, 3 and 4 in \cite{MA}.

\begin{remark} In Examples 2, 3 and 4 of \cite{MA}, it was also imposed the initial and terminal points
constraints: $q_1(t_0), q_2(t_0),q_1(t_1), q_2(t_1)$ are fixed. It
follows from \eqref{ma1} and \eqref{ma2} that, only $q_1(t_0)$ and
$q_2(t_0)$ fixed can guarantee $q_1(t_1)$ and $q_2(t_1)$ also fixed.
\end{remark}
L. Moreau and D. Aeyels \cite{MA} considered the following minimum
problem:

\noindent{\bf (MP)}:\quad To minimize \eqref{ma} subject to
\eqref{ma1} and \eqref{ma2}.

Through defining the augmented state variables
\begin{eqnarray*}
\begin{array}{ll}
&x_1=q_1,\quad x_2=q_2,\quad x_3(t)=\int_{t_0}^ti_3(s)ds,\\
&x_4(t)=\int_{t_0}^ti_5(s)ds,\quad x_5(t)=\int_{t_0}^ti_6(s)ds,
\end{array}
\end{eqnarray*}
\cite{MA} reformulated {\bf (MP)} as an optimal control problem with
terminal state constraints $x_3(t_1)=\lambda_3$,
$x_4(t_1)=\lambda_5$ and $x_5(t_1)=\lambda_6$, and then applying
Pontryagin maximum principle to obtain both the generalized
Hamiltonian and the generalized Euler-Lagrange model of this LC
circuit.

Let us define the symmetric matrix
\begin{eqnarray}
S_1:=\left(\begin{array}{ccc}
L_4+L_3-K_1&-L_4&-L_4\\
-L_4&L_4+L_5-2K_1&L_4\\
-L_4&L_4&L_4+L_6-2K_1
\end{array}\right),
\end{eqnarray}
\begin{eqnarray}
S_2:=\left(\begin{array}{ccc}
L_4+L_3-K_2&-L_4&-L_4\\
-L_4&L_4+L_5-2K_2&L_4\\
-L_4&L_4&L_4+L_6-2K_2
\end{array}\right),
\end{eqnarray}
where $K_1:=\max\{K(C_1),\frac{\widetilde{K}(C_2)}{2}\}$ and
$K_2:=\min\{K(C_1),\frac{\widetilde{K}(C_2)}{2}\}$. Both $K(C_1)$
and $\widetilde{K}(C_2)$ are the unique solutions to
\begin{eqnarray}
\sum_{n=1}^{+\infty}\frac{1}{\frac{2\pi^2C_1}{(t_1-t_0)^2}K(C_1)n^2-\frac{1}{2}}=1,\qquad
K(C_1)>\frac{3(t_1-t_0)^2}{4\pi^2C_1},
\end{eqnarray}
and
\begin{eqnarray}
\sum_{n=1}^{+\infty}\frac{1}{\frac{\pi^2C_2}{(t_1-t_0)^2}\widetilde{K}(C_2)n^2-\frac{1}{2}}=1,\qquad
\widetilde{K}(C_2)>\frac{3(t_1-t_0)^2}{2\pi^2C_2},
\end{eqnarray}
respectively.

In the Appendix, we prove the following
\begin{proposition} It holds that

\begin{description}
\item [(I)]
If the matrix $S_1$ is positively definite, then the Lagrange
functional \eqref{ma} subject to \eqref{ma1} and \eqref{ma2} has a
minimum value at the unique critical point $(i_3^*,i_5^*,i_6^*)\in
C([t_0,t_1],\mathbb{R}^3)$;
\item [(II)]
If the matrix $S_2$ is negatively definite, then the Lagrange
functional \eqref{ma} subject to \eqref{ma1} and \eqref{ma2} has
neither minimum value nor maximum value.
\item [(III)] Let $C_1=C_2$. Then $\widetilde{K}(C_2)=2K(C_1)$, $K_1=K_2$ and $S_1=S_2$.
In these cases, the Lagrange functional \eqref{ma} subject to
\eqref{ma1} and \eqref{ma2} has neither minimum value nor maximum
value provided that the matrix $S_2$ has at least one negative
characteristic root.
\end{description}
\end{proposition}

Define
\begin{eqnarray}\label{matrix}
M:=\left(\begin{array}{ccc}
L_4+L_3&-L_4&-L_4\\
-L_4&L_4+L_5&L_4\\
-L_4&L_4&L_4+L_6
\end{array}\right),
\end{eqnarray}
\begin{eqnarray}\label{matrix1}
N:=\left(\begin{array}{ccc}
\frac{1}{C_1}&0&0\\
0&\frac{1}{C_2}&\frac{1}{C_2}\\
0&\frac{1}{C_2}&\frac{1}{C_2}
\end{array}\right).
\end{eqnarray}
The symmetric matrix $M$ is positively definite due to the
positivity of $L_3,L_4,L_5$ and $L_6$. Hence, we can define the
positively definite matrice $M^{\frac{1}{2}}$ and $M^{-\frac{1}{2}}$
uniquely such that
\begin{eqnarray}\label{matrix2}
M^{\frac{1}{2}}M^{\frac{1}{2}}=M,\qquad
M^{-\frac{1}{2}}M^{-\frac{1}{2}}=M^{-1}.
\end{eqnarray}
There exists an orthogonal matrix $P$ such that
\begin{eqnarray}\label{matrix3}
M^{-\frac{1}{2}}NM^{-\frac{1}{2}}=P^T\left(\begin{array}{ccc}
h_1&0&0\\
0&h_2&0\\
0&0&0
\end{array}\right)P,
\end{eqnarray}
where both $h_1>0$ and $h_2>0$ are the characteristic roots of
$M^{-\frac{1}{2}}NM^{-\frac{1}{2}}$.

\begin{proposition}
It holds that
\begin{description}
\item [(I)]
The Lagrange functional \eqref{ma} subject to \eqref{ma1} and
\eqref{ma2} has a unique critical point $(i_3^*,i_5^*,i_6^*)\in
C([t_0,t_1],\mathbb{R}^3)$ for any
$\lambda_3,\lambda_5,\lambda_6\in\mathbb{R}$ if and only if
\begin{eqnarray}
(t_1-t_0)\sqrt{h_1}\neq k\pi,\text{ and } (t_1-t_0)\sqrt{h_2}\neq
k\pi,\qquad\forall k\in\mathbb{N}^+;
\end{eqnarray}
\item [(II)]
If $(t_1-t_0)\sqrt{h_1}=k\pi$ for some $k\in\mathbb{N}^+$, or
$(t_1-t_0)\sqrt{h_2}=k\pi$ for some $k\in\mathbb{N}^+$, then the
Lagrange functional \eqref{ma} subject to \eqref{ma1} and
\eqref{ma2} has a critical point in $C([t_0,t_1],\mathbb{R}^3)$ for
some $\lambda_3,\lambda_5,\lambda_6\in\mathbb{R}$ if and only if
there exists $\mathbf{b}\in\mathbb{R}^3$ such that
\begin{equation}
\Phi(t_1-t_0)\mathbf{b}=(\lambda_3,\lambda_5,\lambda_6)^T-\int_{t_0}^{t_1}\Phi(t_1-t)M^{-1}\mathbf{a}\,dt,
\end{equation}
where
\begin{equation}
\Phi(t)=M^{-\frac{1}{2}}P^T\left(\begin{array}{ccc}
\frac{1}{\sqrt{h_1}}\sin(\sqrt{h_1}t)&0&0\\
0&\frac{1}{\sqrt{h_2}}\sin(\sqrt{h_2}t)&0\\
0&0&t
\end{array}\right)PM^{\frac{1}{2}},
\end{equation}
\begin{equation}\label{matrix4}
\mathbf{a}:=-(\frac{q_1(t_0)}{C_1}, \frac{q_2(t_0)}{C_2},
\frac{q_2(t_0)}{C_2})^T.
\end{equation}

In these cases, the Lagrange functional \eqref{ma} subject to
\eqref{ma1} and \eqref{ma2} has infinitely many critical points for
these $\lambda_3,\lambda_5,\lambda_6\in\mathbb{R}$.
\end{description}
\end{proposition}
The notion of critical points associated with the Lagrange
functional \eqref{ma} subject to \eqref{ma1} and \eqref{ma2} will be
precisely given by Definition 3.1 and 3.4 in Section 3.

From these propositions, we known the Lagrange functional (2.1)
subject to \eqref{ma1} and \eqref{ma2} has neither minimum nor
maximum value while having critical points in many circumstances.
The critical points in these cases we can understand as
equilibriums. These facts suggest that we should consider the new
variational problem in the next section to replace the minimizing
problem.

\section{A Novel Variational Principle}

In this section, we study the following variational problem:

\begin{eqnarray}\label{y}
J(u):=\int_{t_0}^{t_1}L(x(t),u(t))\,dt=\text{stationary!}
\end{eqnarray}
subject to
\begin{eqnarray}\label{z}
x'=f(x(t),u(t)),\qquad x(t_0)=x_0,
\end{eqnarray}
where $x_0\in\mathbb{R}^n$ is fixed.

It is assumed that
\begin{description}
\item [(AI)] $L:\mathbb{R}^n\times\mathbb{R}^m\mapsto\mathbb{R}$ and
$f:\mathbb{R}^n\times\mathbb{R}^m\mapsto\mathbb{R}^n$ for
$n,m\in\mathbb{N}^+$, are continuously differentiable;
\item [(AII)]
Moreover, for any given $u\in C([t_0,t_1],\mathbb{R}^m)$, the system
\eqref{z} has a unique solution on the whole interval $[t_0,t_1]$,
which will be denoted by $x(\cdot;u)$.
\end{description}

Obviously, the variational equation of \eqref{z} at $(x,u)\in
C([t_0,t_1],\mathbb{R}^n)\times C^1([t_0,t_1],\mathbb{R}^m)$ with
$x=x(\cdot;u)$ is as follows:
\begin{eqnarray}\label{vareq}
\frac{d(\delta x)}{dt}=\frac{\partial f}{\partial
x}(x(t),u(t))\delta x+\frac{\partial f}{\partial u}(x(t),u(t))\delta
u,\qquad \delta x(t_0)=0.
\end{eqnarray}
Let $X(t;u)$ with $t_0\le t\le t_1$ be one fundamental solution
matrix to the homogeneous equation of \eqref{vareq}:
\begin{eqnarray}
\frac{dx}{dt}=\frac{\partial f}{\partial x}(x(t),u(t)) x.
\end{eqnarray}
Define $T(t,s;u):=X(t;u)X^{-1}(s;u)$ with $t_0\le s\le t\le t_1$. By
the variation-of-constants formula, the solution to \eqref{vareq} is
\begin{eqnarray}\label{varcon}
\delta x(t)=\int_{t_0}^tT(t,s;u)\frac{\partial f}{\partial
u}(x(s),u(s))\delta u(s)\,ds,\qquad t_0\le t\le t_1.
\end{eqnarray}

The variation of the functional $J$ of \eqref{y} subject to
\eqref{z} at $u\in C([t_0,t_1];\mathbb{R}^m)$ in the direction $h\in
C([t_0,t_1];\mathbb{R}^m)$ is defined as follows
\begin{eqnarray}\label{u}
\delta J(u;h):=\lim_{\varepsilon\rightarrow 0}\frac{J(u+\varepsilon
h)-J(u)}{\varepsilon},
\end{eqnarray}
in which $J(u)$ is defined by \eqref{y}-\eqref{z}, and
\begin{eqnarray}\label{v}
J(u+\varepsilon h)=\int_{t_0}^{t_1}L(x(t),u(t)+\varepsilon
h(t))\,dt,
\end{eqnarray}
subject to
\begin{eqnarray}\label{w}
x'=f(x(t),u(t)+\varepsilon h(t)),\qquad x(t_0)=x_0.
\end{eqnarray}

Somewhere in the subsequent, for $\Xi=f$, $L$ or $g$, we will denote
$\frac{\partial \Xi}{\partial x}(x(t),u(t))$ (and $\frac{\partial
\Xi}{\partial u}(x(t),u(t))$) simply by $\frac{\partial
\Xi}{\partial x}(t)$ (and $\frac{\partial \Xi}{\partial u}(t)$)
respectively, and analogously for other time variables such as $s$,
$\tau$, etc.

Applying \eqref{varcon} to \eqref{v}-\eqref{w}, we can deduce by
Fubbi Theorem that
\begin{eqnarray}
\begin{array}{ll}
&J(u+\varepsilon h)-J(u)\\
=&\varepsilon\int_{t_0}^{t_1}\frac{\partial L}{\partial
x}(t)\int_{t_0}^tT(t,s;u)\frac{\partial f}{\partial u}(s)h(s)\,ds\,dt+\varepsilon\int_{t_0}^{t_1}\frac{\partial L}{\partial u}(t)h(t)\,dt+o(\varepsilon)\\
=&\varepsilon\int_{t_0}^{t_1}\int_s^{t_1}\frac{\partial L}{\partial
x}(t)T(t,s;u)\,dt\frac{\partial f}{\partial u}(s)h(s)\,ds+\varepsilon\int_{t_0}^{t_1}\frac{\partial L}{\partial u}(t)h(t)\,dt+o(\varepsilon)\\
=&\varepsilon\int_{t_0}^{t_1}\int_t^{t_1}\frac{\partial L}{\partial
x}(s)T(s,t;u)\,ds\frac{\partial f}{\partial u}(t)h(t)\,dt+\varepsilon\int_{t_0}^{t_1}\frac{\partial L}{\partial u}(t)h(t)\,dt+o(\varepsilon)\\
=&\varepsilon\int_{t_0}^{t_1}[\int_t^{t_1}\frac{\partial L}{\partial
x}(s)T(s,t;u)\,ds\frac{\partial f}{\partial u}(t)+\frac{\partial
L}{\partial u}(t)]h(t)\,dt+o(\varepsilon),
\end{array}
\end{eqnarray}
for any $h\in C([t_0,t_1];\mathbb{R}^m)$, and then we have
\begin{eqnarray}\label{proof}
\delta J(u;h)=\int_{t_0}^{t_1}[\int_t^{t_1}\frac{\partial
L}{\partial x}(s)T(s,t;u)\,ds\frac{\partial f}{\partial
u}(t)+\frac{\partial L}{\partial u}(t)]h(t)\,dt.
\end{eqnarray}

\subsection{The first case with additional constraints}

Let us impose some additional constraints
\begin{eqnarray}\label{z1}
\int_{t_0}^{t_1}[Bu(t)+\alpha]\,dt=0,
\end{eqnarray}
 where $B\in\mathbb{R}^{l\times m}$ is a matrix and $\alpha\in\mathbb{R}^l$
is a vector for $l\in\mathbb{N}^+$.

\begin{definition}
\begin{description}
\item [(I)] The admissible set is defined as
\begin{eqnarray}
\mathcal{U}_{ad}=\{u\in
C([t_0,t_1],\mathbb{R}^m)|\int_{t_0}^{t_1}[Bu(t)+\alpha]\,dt=0\};
\end{eqnarray}
\item [(II)] The set of allowed variations is defined as
\begin{eqnarray}
\mathcal{V}_{ad}=\{h\in
C([t_0,t_1],\mathbb{R}^m)|\int_{t_0}^{t_1}Bh(t)\,dt=0\}.
\end{eqnarray}
\end{description}
\end{definition}

\begin{definition} $u\in\mathcal{U}_{ad}$ is called a
critical point for the Lagrange functional \eqref{y} subject to the
equation \eqref{z} and the constraints \eqref{z1} provided that
\begin{eqnarray}\label{critical}
\delta J(u,h)=0,\qquad \forall h\in \mathcal{V}_{ad}.
\end{eqnarray}
In this case, we call the Lagrange functional \eqref{y} subject to
the equation \eqref{z} and the constraints \eqref{z1} is stationary
at this $u\in\mathcal{U}_{ad}$.
\end{definition}

Obviously, $u\in \mathcal{U}_{ad}$ is a critical point for \eqref{y}
subject to \eqref{z} and \eqref{z1} provided that \eqref{y} subject
to \eqref{z} and \eqref{z1} attaches the minimum (or maximum) value
at this $u$.

In classical mechanics, the equation \eqref{z} is the simplest form
as $x'=u$. $\delta u=\delta(x')$ uniquely determine $\delta x$.
Conversely, $\delta x$ also uniquely determine $\delta u$. So we
usually refer the notion of critical points to $x$ instead of $u$ in
that case. For general cases of \eqref{z}, $\delta x$ not always
uniquely determine $\delta u$, which can be discovered from the
dynamic equation \eqref{ma1} of the LC circuit example in Section 2.

The generalization of Hamilton's principle in classical mechanics to
the variational problem of \eqref{y} subject to \eqref{z} and
\eqref{z1} is as follows:

\begin{definition}
$(x,u)$ with $x=x(\cdot,u)$ is called a generalized motion of the
Lagrange functional \eqref{y} subject to the equation \eqref{z} and
the constraints \eqref{z1} provided that $u\in \mathcal{U}_{ad}$ is
a critical point for \eqref{y} subject to \eqref{z} and \eqref{z1}.
\end{definition}
This principle is might as well called {\bf the Hamilton's type
principle}.

\begin{theorem} $(x,u)$ with $x=x(\cdot,u)$ is a generalized motion of the Lagrange functional
\eqref{y} subject to the equation \eqref{z} and the constraints
\eqref{z1}, if and only if $(x,u)$ satisfy both the constraints
\eqref{z1} and the generalized Euler-Lagrange equations
\begin{eqnarray}\label{euler1}
\quad\qquad\begin{cases} \frac{\partial L}{\partial
u}(x(t),u(t))+\int_t^{t_1}\frac{\partial L}{\partial
x}(x(s),u(s))T(s,t;u)\,ds\frac{\partial f}{\partial
u}(x(t),u(t))=\mu^TB,\\
x'(t)=f(x(t),u(t)),\qquad t_0\le t\le t_1,
\end{cases}
\end{eqnarray}
for some $\mu\in\mathbb{R}^l$.
\end{theorem}
{\bf Proof}\quad Let us denote the row vectors of the matrix $B$ by
\begin{eqnarray}
b_i:=(b_{i1},b_{i2},\cdots,b_{im}),\qquad i=1,2,\cdots,l.
\end{eqnarray}
Define that $l$ functions as follows
\begin{eqnarray}
\widehat{b_i}(t)\equiv b_i^T,\qquad t\in [t_0,t_1],
\end{eqnarray}
and
\begin{eqnarray}
L(B):=\text{span}\{\widehat{b_1},\widehat{b_2},\cdots,\widehat{b_l}\},
\end{eqnarray}
which is a complete subspace of $L^2(t_0,t_1;\mathbb{R}^m)$.
$L^2(t_0,t_1;\mathbb{R}^m)=L(B)\oplus L(B)^{\bot}$ where
$L(B)^{\bot}$ is the orthogonal complement space of $L(B)$ in
$L^2(t_0,t_1;\mathbb{R}^m)$. It is well known that
 $C([t_0,t_1],\mathbb{R}^m)$ is
imbedded in $L^2(t_0,t_1;\mathbb{R}^m)$ continuously and densely.
Similarly, $\mathcal{V}_{ad}$ is also imbedded in $L(B)^{\bot}$
continuously and densely. Hence, \eqref{critical} yields that
\eqref{euler1}. \hfill$\Box$\vspace{3mm}

\begin{corollary} $(x,u)$ with $x=x(\cdot;u)$ is a generalized motion of \eqref{y} subject to \eqref{z} if and
only if $(x,u)$ satisfies the generalized Euler-Lagrange equations
\begin{eqnarray}\label{euler}
\begin{cases}
\frac{\partial L}{\partial u}(x(t),u(t))+\int_t^{t_1}\frac{\partial
L}{\partial x}(x(s),u(s))T(s,t;u)\,ds\frac{\partial f}{\partial u}(x(t),u(t))=0,\\
x'(t)=f(x(t),u(t)),\qquad t_0\le t\le t_1.
\end{cases}
\end{eqnarray}
\end{corollary}

\noindent{\bf Example 3.1} In classical mechanics, $f(x,u)=u$ yields
that $\frac{\partial f}{\partial u}(t)\equiv T(s,t;u)\equiv I_n$. If
the initial state $x(t_0)=x_0\in\mathbb{R}^n$ is fixed while the
terminal state $x(t_1)$ is free, then the generalized Euler-Lagrange
equation \eqref{euler} reduces to
\begin{eqnarray}\label{classical}
\frac{\partial L}{\partial
x'}(x(t),x'(t))+\int_t^{t_1}\frac{\partial L}{\partial
x}(x(s),x'(s))\,ds=0,
\end{eqnarray}
as one necessary and sufficient condition for the motion $x\in
C^1([t_0,t_1],\mathbb{R}^n)$.

If the terminal state $x(t_1)=x_1\in\mathbb{R}^n$ is also fixed,
which can be reformulated as the constraints \eqref{z1} with $B=I_n$
and $\alpha=-\frac{x_1-x_0}{t_1-t_0}$. The trajectory $x$ with
$x(t_0)=x_0$ and $x(t_1)=x_1$ is one motion, if and only if $x$
satisfy the equation \eqref{euler1}, which reduces to
\begin{eqnarray}\label{classical2}
\frac{\partial L}{\partial
x'}(x(t),x'(t))+\int_t^{t_1}\frac{\partial L}{\partial
x}(x(s),x'(s))\,ds=\mu,
\end{eqnarray}
for some $\mu\in\mathbb{R}^n$.

Differentiating with respect to $t$, both \eqref{classical} and
\eqref{classical2} yields the classical one
\begin{eqnarray}\label{classical1}
\frac{d}{dt}[\frac{\partial L}{\partial
x'}(x(t),x'(t))]-\frac{\partial L}{\partial x}(x(t),x'(t))=0.
\end{eqnarray}

If some components of the terminal state are fixed while the others
are free, then the generalized Euler-Lagrange equation
\eqref{euler1} is better than the Euler-Lagrange equation
\eqref{classical1} just as the case without terminal state
constraints.

\noindent{\bf Example 3.2} For the LC example of
\eqref{ma}-\eqref{ma1}, let $x:=(q_1,q_2)^T$, $u:=(i_3,i_5,i_6)^T$,
$L(x,u)=\frac{1}{2}L_3i_3^2+\frac{1}{2}L_4(i_3-i_5-i_6)^2+\frac{1}{2}L_5i_5^2+\frac{1}{2}L_6i_6^2-\frac{1}{2C_1}q_1^2-\frac{1}{2C_2}q_2^2$
and
$$
f(x,u)=A\left(\begin{array}{ccc}
i_3\\
i_5\\
i_6
\end{array}\right):=\left(\begin{array}{ccc}
1&0&0\\
0&1&1
\end{array}\right)\left(\begin{array}{ccc}
i_3\\
i_5\\
i_6
\end{array}\right).
$$
Hence $T(t,s;u)\equiv I_2$. The equation \eqref{euler} is
\begin{eqnarray}\label{LC}
\begin{cases}
(L_3+L_4)i_3(t)-L_4i_5(t)-L_4i_6(t)-\int_t^{t_1}\frac{q_1(s)}{C_1}\,ds=0,\\
-L_4i_3(t)+(L_4+L_5)i_5(t)+L_4i_6(t)-\int_t^{t_1}\frac{q_2(s)}{C_2}\,ds=0,\\
-L_4i_3(t)+L_4i_5(t)+(L_4+L_5)i_6(t)-\int_t^{t_1}\frac{q_2(s)}{C_2}\,ds=0;\\
q_1'=i_3,\\
q_2'=i_5+i_6.
\end{cases}
\end{eqnarray}

For the LC example of \eqref{ma}-\eqref{ma1} with the terminal state
$(q_1(t_1),q_2(t_1))$ fixed, can be reformulated as the constraints
\eqref{z1} with $B=A$, the equation \eqref{euler1} is
\begin{eqnarray}\label{LC2}
\qquad\begin{cases}
(L_3+L_4)i_3(t)-L_4i_5(t)-L_4i_6(t)-\int_t^{t_1}\frac{q_1(s)}{C_1}\,ds=\mu_1,\\
-L_4i_3(t)+(L_4+L_5)i_5(t)+L_4i_6(t)-\int_t^{t_1}\frac{q_1(s)}{C_1}\,ds-\int_t^{t_1}\frac{q_2(s)}{C_2}\,ds=\mu_2,\\
-L_4i_3(t)+L_4i_5(t)+(L_4+L_5)i_6(t)-\int_t^{t_1}\frac{q_1(s)}{C_1}\,ds-\int_t^{t_1}\frac{q_2(s)}{C_2}\,ds=\mu_2;\\
q_1'=i_3,\\
q_2'=i_5+i_6,
\end{cases}
\end{eqnarray}
for some $\mu_1,\mu_2\in\mathbb{R}$.

For the LC example of \eqref{ma}-\eqref{ma2}, the constraint
\eqref{ma2} can be reformulated as the constraints \eqref{z1} with
$B=I_3$. The equation \eqref{euler1} is
\begin{eqnarray}\label{LC3}
\qquad\begin{cases}
(L_3+L_4)i_3(t)-L_4i_5(t)-L_4i_6(t)-\int_t^{t_1}\frac{q_1(s)}{C_1}\,ds=\mu_1,\\
-L_4i_3(t)+(L_4+L_5)i_5(t)+L_4i_6(t)-\int_t^{t_1}\frac{q_1(s)}{C_1}\,ds-\int_t^{t_1}\frac{q_2(s)}{C_2}\,ds=\mu_2,\\
-L_4i_3(t)+L_4i_5(t)+(L_4+L_5)i_6(t)-\int_t^{t_1}\frac{q_1(s)}{C_1}\,ds-\int_t^{t_1}\frac{q_2(s)}{C_2}\,ds=\mu_3;\\
q_1'=i_3,\\
q_2'=i_5+i_6,
\end{cases}
\end{eqnarray}
for some $\mu_1,\mu_2,\mu_3\in\mathbb{R}$.

Since the matrix $M$ defined by \eqref{matrix} is positively
definite, by letting $i_4=i_3-i_5-i_6$, all \eqref{LC}, \eqref{LC2}
and \eqref{LC3} yields the same equations
\begin{eqnarray}\label{LC1}
\begin{cases}
L_3i_3'(t)+L_4i_4'(t)+\frac{q_1(t)}{C_1}=0,\\
-L_4i_4'(t)+L_5i_5'(t)+\frac{q_2(t)}{C_2}=0,\\
-L_4i_4'(t)+L_6i_6(t)+\frac{q_2(t)}{C_2}=0;\\
q_1'=i_3,\\
q_2'=i_5+i_6\\
i_4=i_3-i_5-i_6,
\end{cases}
\end{eqnarray}
which is just the generalized Euler-Lagrange equation (3.40)-(3.41)
in \cite{MA}.

\noindent{\bf Example 3.3} For the electromechanical system in
\cite{MA} (see Example 7),
$$
J(u)=\int_{t_0}^{t_1}\{\frac{1}{2}L_1i_1^2+\frac{1}{2}L_2i_2^2+\frac{1}{2}ml^2\omega^2-\frac{q^2}{2C(\theta)}+mgl\cos(\theta)\}\,dt
$$
subject to
$$
\dot{q}=i_1+i_2,\qquad\dot{\theta}=\omega.
$$
where $x:=(q,\theta)^T$, $u:=(i_1,i_2,\omega)^T$. Then
$L(x,u)=\frac{1}{2}L_1i_1^2+\frac{1}{2}L_2i_2^2+\frac{1}{2}ml^2\omega^2-\frac{q^2}{2C(\theta)}+mgl\cos(\theta)$
and
$$
f(x,u)=A\left(\begin{array}{ccc}
i_1\\
i_2\\
\omega
\end{array}\right):=\left(\begin{array}{ccc}
1&1&0\\
0&0&1
\end{array}\right)\left(\begin{array}{ccc}
i_1\\
i_2\\
\omega
\end{array}\right).
$$

If without additional constraints, then the equation \eqref{euler}
is
\begin{eqnarray}\label{LM}
\begin{cases}
L_1i_1(t)-\int_t^{t_1}\frac{q(s)}{C(\theta)}\,ds=0,\\
L_2i_2(t)-\int_t^{t_1}\frac{q(s)}{C(\theta)}\,ds=0,\\
ml^2\omega(t)-\int_t^{t_1}[mgl\sin(\theta)-\frac{q^2(s)}{2C^2(\theta)}C'(\theta)]\,ds=0;\\
q'=i_1+i_2,\\
\theta'=\omega.
\end{cases}
\end{eqnarray}

If the terminal point is fixed, which can be reformulated as the
constraints \eqref{z1} with $B=A$, then the equation \eqref{euler1}
is
\begin{eqnarray}\label{LM1}
\begin{cases}
L_1i_1(t)-\int_t^{t_1}\frac{q(s)}{C(\theta)}\,ds=\mu_1,\\
L_2i_2(t)-\int_t^{t_1}\frac{q(s)}{C(\theta)}\,ds=\mu_1,\\
ml^2\omega(t)-\int_t^{t_1}[mgl\sin(\theta)-\frac{q^2(s)}{2C^2(\theta)}C'(\theta)]\,ds=\mu_2;\\
q'=i_1+i_2,\\
\theta'=\omega,
\end{cases}
\end{eqnarray}
for some $\mu_1,\mu_2\in\mathbb{R}$.

If with the integral constraints
\begin{eqnarray}\label{ma3}
\int_{t_0}^{t_1}i_1\,dt=\lambda_1,\quad\int_{t_0}^{t_1}i_2\,dt=\lambda_2,\quad\int_{t_0}^{t_1}\omega\,dt=\lambda_3,
\end{eqnarray}
which can be reformulated as the constraints \eqref{z1} with
$B=I_3$, and can guarantee the terminal point fixed (similar to
Remark 2.1), then the equation \eqref{euler1} is
\begin{eqnarray}\label{LM2}
\begin{cases}
L_1i_1(t)-\int_t^{t_1}\frac{q(s)}{C(\theta)}\,ds=\mu_1,\\
L_2i_2(t)-\int_t^{t_1}\frac{q(s)}{C(\theta)}\,ds=\mu_2,\\
ml^2\omega(t)-\int_t^{t_1}[mgl\sin(\theta)-\frac{q^2(s)}{2C^2(\theta)}C'(\theta)]\,ds=\mu_3;\\
q'=i_1+i_2,\\
\theta'=\omega,
\end{cases}
\end{eqnarray}
for some $\mu_1,\mu_2,\mu_3\in\mathbb{R}$.

Differentiating the three equations \eqref{LM} , \eqref{LM1} and
\eqref{LM2} yields the same generalized Euler-Lagrange equation
(5.7) in \cite{MA}.

\subsection{The second case with additional constraints}

Consider the Lagrange functional \eqref{y} subject to the dynamic
equation \eqref{z} and some additional constraints
\begin{eqnarray}\label{z2}
\int_{t_0}^{t_1}g(x(t),u(t))\,dt=0,
\end{eqnarray}
where $g:\mathbb{R}^n\times \mathbb{R}^m\mapsto\mathbb{R}^l$ with
$l\in\mathbb{N}^+$ is continuously differentiable.

Let $C([t_0,t_1],\mathbb{R}^m)$ be the Banach space equipped with
the usual maximum norm $\|\cdot\|_C$. Denoted by $v_k\rightharpoonup
v$ in $L^2(t_0,t_1;\mathbb{R}^m)$ the weak convergence of $v_k$ to
$v$ in $L^2(t_0,t_1;\mathbb{R}^m)$.
\begin{definition} For arbitrary given $x_0\in\mathbb{R}^n$,
the admissible set at $x_0$ is defined as
\begin{eqnarray}
\begin{array}{ll}
\mathcal{U}_{ad}(x_0):=&\{u\in C([t_0,t_1],\mathbb{R}^m)| \text{ the
solution } x(\cdot,u) \text{ to } \eqref{z}\\
&\text{ together with } u \text{ satisfy the constraint }
\eqref{z2}\}.
\end{array}
\end{eqnarray}
\end{definition}

Let $B_1$ be the unit ball in $C([t_0,t_1],\mathbb{R}^m)$, which is
a convex set closed under the strong topology of
$L^2(t_0,t_1;\mathbb{R}^m)$. It follows from Mazur theorem that
$B_1$ is also closed under the weak topology in
$L^2(t_0,t_1;\mathbb{R}^m)$. Meanwhile, $B_1$ is weakly precompact
in $L^2(t_0,t_1;\mathbb{R}^m)$.

Hence, for any sequence
$\{u_k\}_{k=1}^{+\infty}\subset\mathcal{U}_{ad}(x_0)$,
$\{\frac{u_k-u}{\|u_k-u\|_C}\}_{k=1}^{+\infty}\subset B_1$ admits a
subsequence weakly convergent to some $h\in B_1$. In this way, we
can define the notion of allowed variation as follows:

\begin{definition}
For any given $x_0\in\mathbb{R}^n$ and $u\in\mathcal{U}_{ad}(x_0)$,
$h\in C([t_0,t_1],\mathbb{R}^m)$ is called an allowed variation
along $u$ at $x_0$ provided that there exists
$\{u_k\}_{k=1}^{+\infty}\subset\mathcal{U}_{ad}(x_0)$ such that
\begin{eqnarray}\label{allow}
\begin{cases}
u_k\rightarrow u \\
\frac{u_k-u}{\|u_k-u\|_C}\rightharpoonup h
\end{cases}
\qquad
\begin{array}{ll}
&\text{in}\quad C([t_0,t_1],\mathbb{R}^m);\\&\text{in}\quad
L^2(t_0,t_1;\mathbb{R}^m).
\end{array}
\end{eqnarray}
The set of all allowed variations along $u$ at $x_0$ is denoted by
$\mathcal{V}_{ad}(x_0,u)$.
\end{definition}

\begin{proposition} For any given $x_0\in\mathbb{R}^n$ and $u\in\mathcal{U}_{ad}(x_0)$,
\begin{eqnarray}\label{con}
\qquad \int_{t_0}^{t_1}[\int_t^{t_1}\frac{\partial g}{\partial
x}(s)T(s,t;u)\,ds\frac{\partial f}{\partial u}(t)+\frac{\partial
g}{\partial u}(t)]h(t)\,dt=0,\qquad\forall h\in
\mathcal{V}_{ad}(x_0,u).
\end{eqnarray}
\end{proposition}

{\bf Proof}\quad

Let $\{u_k\}_{k=1}^{+\infty}\subset\mathcal{U}_{ad}(x_0)$ satisfy
\eqref{allow}, and define $\varepsilon_k:=\|u_k-u\|_C$,
$h_k:=\frac{u_k-u}{\|u_k-u\|_C}$, then

\begin{eqnarray}\label{allow1}
\qquad
\begin{array}{rl}
0=&\int_{t_0}^{t_1}g(x,u_k)\,dt-\int_{t_0}^{t_1}g(x,u)\,dt\\
=&\varepsilon_k\int_{t_0}^{t_1}\frac{\partial g}{\partial
x}(t)\int_{t_0}^tT(t,s;u)\frac{\partial f}{\partial u}(s)h_k(s)\,ds\,dt+\varepsilon\int_{t_0}^{t_1}\frac{\partial g}{\partial u}(t)h_k(t)\,dt+o(\varepsilon_k)\\
=&\varepsilon_k\int_{t_0}^{t_1}\int_s^{t_1}\frac{\partial
g}{\partial
x}(t)T(t,s;u)\,dt\frac{\partial f}{\partial u}(s)h_k(s)\,ds+\varepsilon\int_{t_0}^{t_1}\frac{\partial g}{\partial u}(t)h_k(t)\,dt+o(\varepsilon_k)\\
=&\varepsilon_k\int_{t_0}^{t_1}\int_t^{t_1}\frac{\partial
g}{\partial
x}(s)T(s,t;u)\,ds\frac{\partial f}{\partial u}(t)h_k(t)\,dt+\varepsilon\int_{t_0}^{t_1}\frac{\partial g}{\partial u}(t)h_k(t)\,dt+o(\varepsilon_k)\\
=&\varepsilon_k\int_{t_0}^{t_1}[\int_t^{t_1}\frac{\partial
g}{\partial x}(s)T(s,t;u)\,ds\frac{\partial f}{\partial
u}(t)+\frac{\partial g}{\partial
u}(t)]h_k(t)\,dt+o(\varepsilon_k)\\
=&\varepsilon_k\int_{t_0}^{t_1}[\int_t^{t_1}\frac{\partial
g}{\partial x}(s)T(s,t;u)\,ds\frac{\partial f}{\partial
u}(t)+\frac{\partial g}{\partial u}(t)]h(t)\,dt+o(\varepsilon_k).
\end{array}
\end{eqnarray}
Hence we have \eqref{con}. \hfill$\Box$\vspace{3mm}

Analogous to \eqref{allow1}, let
$\{u_k\}_{k=1}^{+\infty}\subset\mathcal{U}_{ad}(x_0)$ satisfy
\eqref{allow}, we have
\begin{eqnarray}
\qquad
\begin{array}{rl}
J(u_k)-J(u)=&\int_{t_0}^{t_1}L(x,u_k)\,dt-\int_{t_0}^{t_1}L(x,u)\,dt\\
=&\varepsilon_k\int_{t_0}^{t_1}[\int_t^{t_1}\frac{\partial
L}{\partial x}(s)T(s,t;u)\,ds\frac{\partial f}{\partial
u}(t)+\frac{\partial L}{\partial u}(t)]h(t)\,dt+o(\varepsilon_k),
\end{array}
\end{eqnarray}
for any $h\in \mathcal{V}_{ad}(x_0,u)$.

Thus, we define that

\begin{definition} $u\in\mathcal{U}_{ad}(x_0)$ is called a
critical point for the Lagrange functional \eqref{y} subject to the
equation \eqref{z} and the constraints \eqref{z2} provided that
\begin{eqnarray}\label{aa}
\delta J(u;h)=0,\qquad \forall h\in \mathcal{V}_{ad}(x_0,u).
\end{eqnarray}
In this case, we call \eqref{y} subject to \eqref{z} and \eqref{z2}
is stationary at this $u\in \mathcal{U}_{ad}(x_0)$.
\end{definition}

\begin{definition}
$(x,u)$ with $x=x(\cdot,u)$ is called a generalized motion of the
Lagrange functional \eqref{y} subject to the equation \eqref{z} and
the constraints \eqref{z2} provided that $u\in
\mathcal{U}_{ad}(x_0)$ is a critical point for \eqref{y} subject to
\eqref{z} and \eqref{z2}.
\end{definition}

\begin{theorem} $(x,u)$ with $x=x(\cdot,u)$ is a generalized motion of the Lagrange functional
\eqref{y} subject to the equation \eqref{z} and the constraints
\eqref{z2} provided that $(x,u)$ satisfy \eqref{z2} and the
generalized Euler-Lagrange equations
\begin{eqnarray}\label{euler2}
\begin{cases}\frac{\partial L}{\partial
u}(t)-\mu^T\frac{\partial g}{\partial
u}(t)+\int_t^{t_1}[\frac{\partial L}{\partial
x}(s)-\mu^T\frac{\partial g}{\partial
x}(s)]T(s,t;u)\,ds\frac{\partial f}{\partial
u}(t)=0,\\
x'(t)=f(x(t),u(t)),\qquad t_0\le t\le t_1,
\end{cases}
\end{eqnarray}
for some $\mu\in\mathbb{R}^l$.
\end{theorem}
{\bf Proof}\quad Combining the generalized Euler-Lagrange equation
\eqref{euler2} and \eqref{con} yields \eqref{aa}. The proof is
completed. \hfill$\Box$\vspace{3mm}

\section{The canonical Hamiltonian systems}

\subsection{The case without additional
constraints}

\noindent{\bf (AIIIa)} The equation
\begin{eqnarray}
0=p^T\frac{\partial f}{\partial u}(x,u)-\frac{\partial L}{\partial
u}(x,u),
\end{eqnarray}
admits one smooth solution
\begin{eqnarray}\label{leg}
u=\varphi(x,p),\qquad\forall (x,p)\in\mathbb{R}^n\times\mathbb{R}^n.
\end{eqnarray}

Define the pseudo Hamiltonian
\begin{eqnarray}\label{pham}
\mathscr{H}(x,p,u):=p^Tf(x,u)-L(x,u),\qquad\forall
(x,p,u)\in\mathbb{R}^n\times\mathbb{R}^n\times\mathbb{R}^m,
\end{eqnarray}
and define under the assumption {\bf (AIIIa)} the Hamiltonian
\begin{eqnarray}\label{ham}
H(x,p):=\mathscr{H}(x,p,\varphi(x,p)),\qquad\forall
(x,p)\in\mathbb{R}^n\times\mathbb{R}^n,
\end{eqnarray}
which is called the generalized Legendre transformation of
$\mathscr{H}$.

\begin{remark} The Legendre transformation of
$\mathscr{H}$ is defined as
$$
H(x,p):=\max_{u\in\mathbb{R}^m}\mathscr{H}(x,p,u),\qquad\forall
(x,p)\in\mathbb{R}^n\times\mathbb{R}^n.
$$
The definition of \eqref{ham} is adapted from {\bf (A1)} in
\cite{SW} (p.39), which indicates there maybe exist different
Hamiltonian descriptions for the same equilibrium in nonlinear LC
circuits.
\end{remark}

Suppose that $(x,u)$ is one solution to the generalized
Euler-Lagrange equation \eqref{euler}.

Let
\begin{eqnarray}\label{p}
p(t)^T:=-\int_t^{t_1}\frac{\partial L}{\partial
x}(x(s),u(s))T(s,t;u)\,ds,\qquad t\in [t_0,t_1],
\end{eqnarray}
then \eqref{euler} can be recast as
\begin{eqnarray}\label{re}
\begin{cases}
p(t)^T\frac{\partial f}{\partial u}(x(t),u(t))-\frac{\partial L}{\partial u}(x(t),u(t))=0,\\
x'(t)=f(x(t),u(t)),\qquad t\in [t_0,t_1].
\end{cases}
\end{eqnarray}

It follows from the first equality of \eqref{re} and the definitions
of $H$ and $\mathscr{H}$ that
\begin{eqnarray}\label{h}
\begin{cases}
\frac{\partial\mathscr{H}}{\partial
p}(x(t),p(t),u(t))\equiv\frac{\partial H}{\partial p}(x(t),p(t)),\\
\frac{\partial\mathscr{H}}{\partial
x}(x(t),p(t),u(t))\equiv\frac{\partial H}{\partial x}(x(t),p(t)),
\end{cases}
\qquad t\in [t_0,t_1].
\end{eqnarray}
Meanwhile, the second equality of \eqref{re} yields that
\begin{eqnarray}\label{h1}
x'(t)=f(x(t),u(t))=\frac{\partial\mathscr{H}}{\partial
p}(x(t),p(t),u(t)),
\end{eqnarray}
and differentiating \eqref{p} yields that
\begin{eqnarray}\label{h2}
\begin{array}{ll}
[p'(t)]^T&=\frac{\partial L}{\partial
x}(x(t),u(t))+\int_t^{t_1}\frac{\partial L}{\partial
x}(x(s),u(s))T(s,t;u)\,ds\frac{\partial f}{\partial
x}(x(t),u(t))\\
&=-\frac{\partial\mathscr{H}}{\partial x}(x(t),p(t),u(t)).
\end{array}
\end{eqnarray}

Hence, we have
\begin{theorem} Let the assumption {\bf (AIIIa)} holds. Suppose that $(x,u)$ is one solution to the generalized Euler-Lagrange equation
\eqref{euler}, then $(x,p)$ given by \eqref{p} is one solution to
the canonical Hamiltonian system
\begin{eqnarray}\label{hamsys}
\begin{cases}
x'(t)=\nabla_pH(x,p),\\
p'(t)=-\nabla_xH(x,p),
\end{cases}
\begin{array}{ll}
&x(t_0)=x_0,\\
&p(t_1)=0.
\end{array}
\end{eqnarray}
Conversely, if $(x,p)$ is one solution to \eqref{hamsys}, then
$(x,u)$ given by \eqref{leg} is also one solution to the generalized
Euler-Lagrange equation \eqref{euler}.
\end{theorem}

\subsection{The cases with additional special constraints}
In \eqref{z2}, let us assume that there exists some matrix
$Q\in\mathbb{R}^{l\times n}$ and $\beta\in\mathbb{R}^l$ such that
\begin{eqnarray}\label{assum}
g(x,u)=Qf(x,u)+\beta,\qquad\forall
(x,u)\in\mathbb{R}^n\times\mathbb{R}^m.
\end{eqnarray}

\begin{remark} The terminal state constraints with $x(t_1)=x_1$
fixed, can be reformulated as \eqref{assum} with $Q=I_n$ and
$\beta=-\frac{x_1-x_0}{t_1-t_0}$.

If there exist $A\in\mathbb{R}^{n\times m}$,
$B\in\mathbb{R}^{l\times m}$ and $\alpha\in\mathbb{R}^l$ such that
$f(x,u)=Au$ and $g(x,u)=Bu+\alpha$, then \eqref{assum} is equivalent
to
\begin{eqnarray}
\text{\rm rank}(A)=\text{\rm rank}\left(\begin{array}{ccc}
A\\
B
\end{array}\right).
\end{eqnarray}
\end{remark}

Let
\begin{eqnarray}\label{p2}
p(t)^T:=-\int_t^{t_1}[\frac{\partial L}{\partial
x}(s)-\mu^T\frac{\partial g}{\partial x}(s)]T(s,t;u)\,ds+\mu^T Q,
\end{eqnarray}
for the parameters $\mu=(\mu_1,\cdots\mu_l)^T\in\mathbb{R}^l$. Then,
\eqref{euler2} can be recast as
\begin{eqnarray}\label{re2}
\begin{cases}
p(t)^T\frac{\partial f}{\partial u}(x(t),u(t))-\frac{\partial L}{\partial u}(x(t),u(t))=0,\\
x'(t)=f(x(t),u(t)).
\end{cases}
\end{eqnarray}

\begin{theorem} Assume that both \eqref{assum} and the assumption {\bf (AIIIa)} holds.
Suppose that $(x,u)$ is one solution to the generalized
Euler-Lagrange equations \eqref{euler2}, then $(x,p)$ given by
\eqref{p2} is one solution to the Hamiltonian system
\begin{eqnarray}\label{hamsys2}
\begin{cases}
x'(t)=\nabla_pH(x,p),\\
p'(t)=-\nabla_xH(x,p),
\end{cases}
\begin{array}{ll}
&x(t_0)=x_0,\\
&p(t_1)=Q^T\mu,
\end{array}
\end{eqnarray}
where the Hamiltonian $H$ is defined in \eqref{ham}.

Conversely, if $(x,p)$ is one solution to \eqref{hamsys2}, then
$(x,u)$ given by \eqref{leg} is also one solution to the generalized
Euler-Lagrange equation \eqref{euler2} for the parameters
$\mu\in\mathbb{R}^l$.
\end{theorem}

\subsection{The general cases}
Define the pseudo Hamiltonian
\begin{eqnarray}\label{pham1}
\begin{array}{ll}
\mathscr{H}(x,p,u;\mu):=&p^Tf(x,u)-L(x,u)+\mu^Tg(x,u),\\
&\qquad\forall
(x,p,u,\mu)\in\mathbb{R}^n\times\mathbb{R}^n\times\mathbb{R}^m\times\mathbb{R}^l.
\end{array}
\end{eqnarray}

\noindent{\bf (AIIIb)} The equation
\begin{eqnarray}
0=p^T\frac{\partial f}{\partial u}(x,u)-\frac{\partial L}{\partial
u}(x,u)+\mu^T\frac{\partial g}{\partial u}(x,u),
\end{eqnarray}
admits one smooth solution
\begin{eqnarray}\label{leg1}
u=\varphi(x,p;\mu),\qquad\forall
(x,p)\in\mathbb{R}^n\times\mathbb{R}^n,
\end{eqnarray}
for some parameters $\mu=(\mu_1,\cdots,\mu_l)^T\in\mathbb{R}^l$.

Under the assumption {\bf (AIIIb)}, let us define the Hamiltonian
\begin{eqnarray}\label{ham1}
H(x,p;\mu):=\mathscr{H}(x,p,\varphi(x,p;\mu);\mu),\qquad\forall
(x,p)\in\mathbb{R}^n\times\mathbb{R}^n.
\end{eqnarray}

Suppose that $(x,u)$ is one solution to the generalized
Euler-Lagrange equation \eqref{euler2} for these parameters
$\mu\in\mathbb{R}^l$. Let
\begin{eqnarray}\label{p1}
p(t)^T:=-\int_t^{t_1}[\frac{\partial L}{\partial
x}(x(s),u(s))-\mu^T\frac{\partial g}{\partial
x}(x(s),u(s))]T(s,t;u)\,ds,
\end{eqnarray}
then the generalized Euler-Lagrange equation \eqref{euler2} can be
recast as
\begin{eqnarray}\label{re1}
\begin{cases}
p(t)^T\frac{\partial f}{\partial u}(x(t),u(t))-\frac{\partial L}{\partial u}(x(t),u(t))+\mu^T\frac{\partial g}{\partial u}(x(t),u(t))=0,\\
x'(t)=f(x(t),u(t)),\qquad t\in [t_0,t_1].
\end{cases}
\end{eqnarray}

Similar to the proof of Theorem 4.1, we have
\begin{theorem} Let the assumptions {\bf (AIIIb)} holds.  Suppose that $(x,u)$ is one solution to the generalized
Euler-Lagrange equation \eqref{euler2} for the parameters
$\mu\in\mathbb{R}^l$, then $(x,p)$ given by \eqref{p1} is one
solution to the Hamiltonian system
\begin{eqnarray}\label{hamsys1}
\begin{cases}
x'(t)=\nabla_pH(x,p;\mu),\\
p'(t)=-\nabla_xH(x,p;\mu),
\end{cases}
\begin{array}{ll}
&x(t_0)=x_0,\\
&p(t_1)=0.
\end{array}
\end{eqnarray}

Conversely, if $(x,p)$ is one solution to \eqref{hamsys1} for the
parameters $\mu\in\mathbb{R}^l$, then $(x,u)$ given by \eqref{leg1}
is also one solution to the generalized Euler-Lagrange equation
\eqref{euler2}.
\end{theorem}

\begin{remark}
From Theorem 4.1, Theorem 4.2 and Theorem 4.3, we can find out that
the canonical Hamiltonian is not enough to describe the energy
function and the dynamic equation when the constraints \eqref{z2}
become more complex.
\end{remark}

\noindent{\bf Example 4.1} (Continued from Example 3.1) The
Hamiltonian is as usual defined by
$H(x,p)=\max_{u\in\mathbb{R}^n}\{p^Tu-L(x,u)\}$, and it follows from
Theorem 4.1, that the dynamic system without additional constraints
is described by the two-point boundary value problem of the
canonical Hamiltonian system.

Two different description of the Hamiltonian system. The first is
the classical approach. Applying Theorem 4.2 yields
\begin{eqnarray}
\begin{cases}
x'(t)=\nabla_pH(x,p),\\
p'(t)=-\nabla_xH(x,p),
\end{cases}
\begin{array}{ll}
&x(t_0)=x_0\\
&p(t_1)=p_1,
\end{array}
\end{eqnarray}
where the terminal costate $p_1\in\mathbb{R}^n$ are the parameters
such that the terminal state constraints $x(t_1)=x_1$ satisfied. The
second is applying Theorem 4.3.
\begin{eqnarray}
\begin{cases}
x'(t)=\nabla_p\widetilde{H}(x,p;\mu),\\
p'(t)=-\nabla_x\widetilde{H}(x,p;\mu),
\end{cases}
\begin{array}{ll}
&x(t_0)=x_0\\
&p(t_1)=0,
\end{array}
\end{eqnarray}
where the Hamiltonian
$\widetilde{H}(x,p;\mu)=\max_{u\in\mathbb{R}^n}\{p^Tu-L(x,u)+\mu^T
u\}$, and $\mu\in\mathbb{R}^n$ are the parameters such that the
terminal state constraints $x(t_1)=x_1$ satisfied.

In fact, $H$ is the energy function and $\widetilde{H}$ is the
energy function with constraints.

\noindent{\bf Example 4.2} (Continued from Example 3.2) The
Hamiltonian is
\begin{eqnarray*}
\begin{array}{ll}
H(x,p)&=\max_{(i_3,i_5,i_6)\in\mathbb{R}^3}\{p_1i_3+p_2(i_5+i_6)-\frac{1}{2}L_3i_3^2\\
&\qquad-\frac{1}{2}L_4(i_3-i_5-i_6)^2-\frac{1}{2}L_5i_5^2-\frac{1}{2}L_6i_6^2+\frac{1}{2C_1}q_1^2+\frac{1}{2C_2}q_2^2\}\\
&=\frac{1}{2}(p_1,p_2,p_2)M^{-1}\left(\begin{array}{ccc}
p_1\\
p_2\\
p_2
\end{array}\right)+\frac{1}{2C_1}q_1^2+\frac{1}{2C_2}q_2^2,
\end{array}
\end{eqnarray*}
where $M$ is the positively definite matrix defined by
\eqref{matrix}. Then applying Theorem 4.1 and Theorem 4.2 to this
Hamiltonian, we obtain the canonical Hamiltonian systems description
of this model without additional constraints and with the terminal
state constraints, respectively.

The Hamiltonian is
\begin{eqnarray*}
\begin{array}{ll}
H(x,p;\mu)&=\max_{(i_3,i_5,i_6)\in\mathbb{R}^3}\{(p_1+\mu_1)i_3+(p_2+\mu_2)i_5+(p_2+\mu_3)i_6\\
&\qquad-\frac{1}{2}L_3i_3^2-\frac{1}{2}L_4(i_3-i_5-i_6)^2-\frac{1}{2}L_5i_5^2-\frac{1}{2}L_6i_6^2+\frac{1}{2C_1}q_1^2+\frac{1}{2C_2}q_2^2\}\\
&=\frac{1}{2}(p_1+\mu_1,p_2+\mu_2,p_2+\mu_3)M^{-1}\left(\begin{array}{ccc}
p_1+\mu_1\\
p_2+\mu_2\\
p_2+\mu_3
\end{array}\right)+\frac{1}{2C_1}q_1^2+\frac{1}{2C_2}q_2^2,
\end{array}
\end{eqnarray*}
where $M$ is the positively definite matrix defined by
\eqref{matrix}. Then applying Theorem 4.3 to this Hamiltonian, we
obtain the canonical Hamiltonian systems description of this model.
Intuitively, the original energy function is not enough in general
to describe the dynamic system with constraints since the
constraints involves three parameters $\mu_1,\mu_2,\mu_3$ while the
dimension of the costate is only 2.

\noindent{\bf Example 4.3} (Continued from Example 3.3) The
Hamiltonian is
\begin{eqnarray*}
\begin{array}{ll}
H(x,p)&=\max_{(i_1,i_2,\omega)\in\mathbb{R}^3}\{p_1(i_1+i_2)+p_2\omega-\frac{1}{2}L_1i_1^2-\frac{1}{2}L_2i_2^2-\frac{1}{2}ml^2\omega^2\\
&\qquad+\frac{q^2}{2C(\theta)}-mgl\cos(\theta)\}\\
&=\frac{1}{2}(\frac{1}{L_1}+\frac{1}{L_2})p_1^2+\frac{1}{2ml^2}p_2^2+\frac{q^2}{2C(\theta)}-mgl\cos(\theta).
\end{array}
\end{eqnarray*}
Then applying Theorem 4.1 and Theorem 4.2 to this Hamiltonian, we
obtain the canonical Hamiltonian systems description of this model
without additional constraints and with the terminal state
constraints, respectively.

The Hamiltonian is
\begin{eqnarray*}
\begin{array}{ll}
H(x,p)&=\max_{(i_1,i_2,\omega)\in\mathbb{R}^3}\{(p_1+\mu_1)i_1+(p_1+\mu_2)i_2+(p_2+\mu_3)\omega\\
&\qquad-\frac{1}{2}L_1i_1^2-\frac{1}{2}L_2i_2^2-\frac{1}{2}ml^2\omega^2+\frac{q^2}{2C(\theta)}-mgl\cos(\theta)\}\\
&=\frac{1}{2L_1}(p_1+\mu_1)^2+\frac{1}{2L_2}(p_1+\mu_2)^2+\frac{1}{2ml^2}(p_2+\mu_3)^2+\frac{q^2}{2C(\theta)}-mgl\cos(\theta).
\end{array}
\end{eqnarray*}
Then applying Theorem 4.3 to this Hamiltonian, we obtain the
canonical Hamiltonian systems description of this model. Similarly,
the original energy function is also not enough in general to
describe the dynamic system with these constraints \eqref{ma3}.

\section{Appendix}

Let $L^2(t_0,t_1;\mathbb{R})$ be the Hilbert space equipped with the
inner product
\begin{eqnarray}
\langle
u,v\rangle:=\frac{2}{t_1-t_0}\int_{t_0}^{t_1}u(t)v(t)dt,\qquad
\forall u,v\in L^2(t_0,t_1;\mathbb{R}),
\end{eqnarray}
and define $e_0$, $e_n$, $\widetilde{e_n}\in
L^2(t_0,t_1;\mathbb{R})$ for $n\in \mathbb{N}^+$ as follows:
$e_0(t)\equiv\frac{1}{\sqrt{2}}$ and
\begin{eqnarray}
e_n(t):=\cos\frac{2n\pi(t-t_0)}{t_1-t_0},\quad
\widetilde{e_n}(t):=\sin\frac{2n\pi(t-t_0)}{t_1-t_0}, \qquad t\in
[t_0,t_1],
\end{eqnarray}
then
$\{e_0,e_1,\widetilde{e_1},e_2,\widetilde{e_2},\cdots,e_n,\widetilde{e_n},\cdots\}$
is an orthonormal basis of $L^2(t_0,t_1;\mathbb{R})$.

The constraint \eqref{ma2} yields that
\begin{eqnarray}
i_k=\frac{\sqrt{2}\lambda_k}{t_1-t_0}e_0+\sum_{n=1}^{+\infty}(a_{k,n}e_n+b_{k,n}\widetilde{e_n}),\quad
\end{eqnarray}
with $\sum_{n=1}^{+\infty}(a_{k,n}^2+b_{k,n}^2)<+\infty$, for
$k=3,5,6$.

Through direct calculations, we have
\begin{lemma} The Lagrange functional \eqref{ma} subject to \eqref{ma1} and \eqref{ma2} can be represented as follows:
\begin{eqnarray}
J(i_3,i_5,i_6)=\mathbf{Q}+\mathbf{L}+\mathbf{N},
\end{eqnarray}
where
\begin{eqnarray}\begin{array}{ll}
\mathbf{Q}=&\quad\{\frac{(t_1-t_0)L_4}{4}\sum_{n=1}^{+\infty}(a_{3,n}-a_{5,n}-a_{6,n})^2\\
&+[\frac{(t_1-t_0)L_3}{4}\sum_{n=1}^{+\infty}a_{3,n}^2-\frac{(t_1-t_0)^3}{16\pi^2C_1}\sum_{n=1}^{+\infty}(\frac{a_{3,n}}{n})^2]\\
&+[\frac{(t_1-t_0)L_5}{4}\sum_{n=1}^{+\infty}a_{5,n}^2+\frac{(t_1-t_0)L_6}{4}\sum_{n=1}^{+\infty}a_{6,n}^2
-\frac{(t_1-t_0)^3}{16\pi^2C_2}\sum_{n=1}^{+\infty}(\frac{a_{5,n}+a_{6,n}}{n})^2]\}\\
&\quad+\{\frac{(t_1-t_0)L_4}{4}\sum_{n=1}^{+\infty}(b_{3,n}-b_{5,n}-b_{6,n})^2\\
&+[\frac{(t_1-t_0)L_3}{4}\sum_{n=1}^{+\infty}b_{3,n}^2-\frac{(t_1-t_0)^3}{16\pi^2C_1}\sum_{n=1}^{+\infty}(\frac{b_{3,n}}{n})^2
-\frac{(t_1-t_0)^3}{8\pi^2C_1}(\sum_{n=1}^{+\infty}\frac{b_{3,n}}{n})^2]\\
&+[\frac{(t_1-t_0)L_5}{4}\sum_{n=1}^{+\infty}b_{5,n}^2
+\frac{(t_1-t_0)L_6}{4}\sum_{n=1}^{+\infty}b_{6,n}^2\\
&\quad-\frac{(t_1-t_0)^3}{16\pi^2C_2}\sum_{n=1}^{+\infty}(\frac{b_{5,n}+b_{6,n}}{n})^2
-\frac{(t_1-t_0)^3}{8\pi^2C_2}(\sum_{n=1}^{+\infty}\frac{b_{5,n}+b_{6,n}}{n})^2]\},\end{array}
\end{eqnarray}
\begin{eqnarray}\begin{array}{ll}
\mathbf{L}=&-\frac{(t_1-t_0)^2}{4\pi C_1}[2q_1(t_0)+\lambda_3]\sum_{n=1}^{+\infty}\frac{b_{3,n}}{n}\\
&-\frac{(t_1-t_0)^2}{4\pi
C_2}[2q_2(t_0)+\lambda_5+\lambda_6]\sum_{n=1}^{+\infty}\frac{b_{5,n}+b_{6,n}}{n},
\end{array}
\end{eqnarray}
\begin{eqnarray}\begin{array}{ll}
\mathbf{N}=&\frac{1}{2(t_1-t_0)}[L_3\lambda_3^2+L_5\lambda_5^2+L_6\lambda_6^2+L_4(\lambda_3-\lambda_5-\lambda_6)^2]\\
&-\frac{t_1-t_0}{6C_1}[3q_1(t_0)^2+3q_1(t_0)\lambda_3+\lambda_3^2]\\
&-\frac{t_1-t_0}{6C_2}[3q_2(t_0)^2+3q_2(t_0)(\lambda_5+\lambda_6)+(\lambda_5+\lambda_6)^2].
\end{array}
\end{eqnarray}
\end{lemma}

\begin{lemma} Let $\alpha>0$ and $\beta>0$. Then
\begin{eqnarray}\label{app}
\alpha\sum_{n=1}^{+\infty}\frac{x_n^2}{n^2}+\beta(\sum_{n=1}^{+\infty}\frac{x_n}{n})^2\le
K\sum_{n=1}^{+\infty}x_n^2,
\end{eqnarray}
where $K$ is the unique solution to the equation
\begin{eqnarray}\label{app1}
\beta \sum_{n=1}^{+\infty}\frac{1}{Kn^2-\alpha}=1,\qquad
K>\alpha+\beta.
\end{eqnarray}

The equality in \eqref{app} holds if and only if
$x_n=\frac{C}{Kn-\frac{\alpha}{n}}$ for arbitrary $C\in \mathbb{R}$.
\end{lemma}
{\bf Proof}
\begin{eqnarray*}
K=\sup_{\sum_{n=1}^{+\infty}x_n^2=1}[\alpha\sum_{n=1}^{+\infty}\frac{x_n^2}{n^2}+\beta(\sum_{n=1}^{+\infty}\frac{x_n}{n})^2]\ge\alpha+\beta.
\end{eqnarray*}

For $l\in\mathbb{N}^+$, it can be shown by Lagrange multiplier rule
that the problem:
\begin{eqnarray*}
\left\{\begin{array}{ll}&\text{To maximize
}\quad \alpha\sum_{n=1}^l\frac{x_n^2}{n^2}+\beta(\sum_{n=1}^l\frac{x_n}{n})^2\\
&\text{subject to }\quad\sum_{n=1}^lx_n^2=1\end{array}\right.
\end{eqnarray*}
has only two solutions
\begin{eqnarray*}
x_{l,n}=\pm[\sum_{n=1}^l\frac{1}{(k_ln-\frac{\alpha}{n})^2}]^{-\frac{1}{2}}\frac{1}{k_ln-\frac{\alpha}{n}},\qquad
n=1,2,\cdots,l,
\end{eqnarray*}
where
\begin{eqnarray*}
k_l=\max_{\sum_{n=1}^lx_n^2=1}[\alpha\sum_{n=1}^l\frac{x_n^2}{n^2}+\beta(\sum_{n=1}^l\frac{x_n}{n})^2],
\end{eqnarray*}
satisfies
\begin{eqnarray*}\label{equation1}
\beta \sum_{n=1}^l\frac{1}{k_ln^2-\alpha}=1.
\end{eqnarray*}

Obviously, it follows from $\lim_{l\rightarrow +\infty}k_l=K$ that,
the equation \eqref{app1} and
\begin{eqnarray*}
x_n^*:=\lim_{l\rightarrow
+\infty}x_{l,n}=\pm[\sum_{n=1}^{+\infty}\frac{1}{(Kn-\frac{\alpha}{n})^2}]^{-\frac{1}{2}}\frac{1}{Kn-\frac{\alpha}{n}},\qquad
n\in \mathbb{N}^+,
\end{eqnarray*}
\begin{eqnarray*}
\alpha\sum_{n=1}^{+\infty}\frac{{x_n^*}^2}{n^2}+\beta(\sum_{n=1}^{+\infty}\frac{x_n^*}{n})^2
=\lim_{l\rightarrow
+\infty}[\alpha\sum_{n=1}^l\frac{x_{l,n}^2}{n^2}+\beta(\sum_{n=1}^l\frac{x_{l,n}}{n})^2]=K,
\end{eqnarray*}
which implies the sufficiency.

The necessity can be shown directly by Ljusternik Theorem (the
Lagrange multiplier rule in infinite dimensional space, see pp.290
in \cite{Ze}). \hfill$\Box$\vspace{3mm}

Similarly, we have
\begin{lemma} Let $\alpha>0$ and $\beta>0$. Then
\begin{eqnarray}\label{app2}
\alpha\sum_{n=1}^{+\infty}\frac{(x_n+y_n)^2}{n^2}
+\beta(\sum_{n=1}^{+\infty}\frac{x_n+y_n}{n})^2\le
\widetilde{K}(\sum_{n=1}^{+\infty}x_n^2+\sum_{n=1}^{+\infty}y_n^2),
\end{eqnarray}
where $\widetilde{K}$ is the unique solution to the equation
\begin{eqnarray}
2\beta
\sum_{n=1}^{+\infty}\frac{1}{\widetilde{K}n^2-2\alpha}=1,\qquad
\widetilde{K}>2(\alpha+\beta).
\end{eqnarray}

The equality in \eqref{app2} holds if and only if
$x_n=y_n=\frac{C}{\widetilde{K}n-\frac{2\alpha}{n}}$ for arbitrary
$C\in \mathbb{R}$. \hfill$\Box$\vspace{3mm}
\end{lemma}

\vskip 3mm

{\bf Proof of Proposition 2.1}\quad  Let
$\mathbf{Q}:=\mathbf{Q}_1+\mathbf{Q}_2$, where
\begin{eqnarray}\begin{array}{ll}
\qquad\quad\mathbf{Q}_1=&\quad\frac{(t_1-t_0)L_4}{4}\sum_{n=1}^{+\infty}(a_{3,n}-a_{5,n}-a_{6,n})^2\\
&+[\frac{(t_1-t_0)L_3}{4}\sum_{n=1}^{+\infty}a_{3,n}^2-\frac{(t_1-t_0)^3}{16\pi^2C_1}\sum_{n=1}^{+\infty}(\frac{a_{3,n}}{n})^2]\\
&+[\frac{(t_1-t_0)L_5}{4}\sum_{n=1}^{+\infty}a_{5,n}^2+\frac{(t_1-t_0)L_6}{4}\sum_{n=1}^{+\infty}a_{6,n}^2
-\frac{(t_1-t_0)^3}{16\pi^2C_2}\sum_{n=1}^{+\infty}(\frac{a_{5,n}+a_{6,n}}{n})^2],\end{array}
\end{eqnarray}
\begin{eqnarray}\begin{array}{ll}
\quad\mathbf{Q}_2=&\quad\frac{(t_1-t_0)L_4}{4}\sum_{n=1}^{+\infty}(b_{3,n}-b_{5,n}-b_{6,n})^2\\
&+[\frac{(t_1-t_0)L_3}{4}\sum_{n=1}^{+\infty}b_{3,n}^2-\frac{(t_1-t_0)^3}{16\pi^2C_1}\sum_{n=1}^{+\infty}(\frac{b_{3,n}}{n})^2
-\frac{(t_1-t_0)^3}{8\pi^2C_1}(\sum_{n=1}^{+\infty}\frac{b_{3,n}}{n})^2]\\
&+[\frac{(t_1-t_0)L_5}{4}\sum_{n=1}^{+\infty}b_{5,n}^2
+\frac{(t_1-t_0)L_6}{4}\sum_{n=1}^{+\infty}b_{6,n}^2\\
&\quad-\frac{(t_1-t_0)^3}{16\pi^2C_2}\sum_{n=1}^{+\infty}(\frac{b_{5,n}+b_{6,n}}{n})^2
-\frac{(t_1-t_0)^3}{8\pi^2C_2}(\sum_{n=1}^{+\infty}\frac{b_{5,n}+b_{6,n}}{n})^2].\end{array}
\end{eqnarray}

{\bf (I)} By Lemma 6.2 and 6.3,
\begin{eqnarray}\label{base}
\begin{array}{rl} \frac{4}{t_1-t_0}\mathbf{Q}_2
=& L_4\sum_{n=1}^{+\infty}(b_{3,n}-b_{5,n}-b_{6,n})^2\\
&+[L_3\sum_{n=1}^{+\infty}b_{3,n}^2-\frac{(t_1-t_0)^2}{4\pi^2C_1}\sum_{n=1}^{+\infty}(\frac{b_{3,n}}{n})^2
-\frac{(t_1-t_0)^2}{2\pi^2C_1}(\sum_{n=1}^{+\infty}\frac{b_{3,n}}{n})^2]\\
&+[L_5\sum_{n=1}^{+\infty}b_{5,n}^2
+L_6\sum_{n=1}^{+\infty}b_{6,n}^2\\
&\quad-\frac{(t_1-t_0)^2}{4\pi^2C_2}\sum_{n=1}^{+\infty}(\frac{b_{5,n}+b_{6,n}}{n})^2
-\frac{(t_1-t_0)^2}{2\pi^2C_2}(\sum_{n=1}^{+\infty}\frac{b_{5,n}+b_{6,n}}{n})^2]\\
\ge& L_4\sum_{n=1}^{+\infty}(b_{3,n}-b_{5,n}-b_{6,n})^2\\
&+[L_3-K_1]\sum_{n=1}^{+\infty}b_{3,n}^2+[L_5-2K_1]\sum_{n=1}^{+\infty}b_{5,n}^2+[L_6-2K_1]\sum_{n=1}^{+\infty}b_{6,n}^2\\
=&\sum_{n=1}^{+\infty}(b_{3,n},b_{5,n},b_{6,n})S_1(b_{3,n},b_{5,n},b_{6,n})^T.\end{array}
\end{eqnarray}
If $S_1$ is positively definite, then it follows from \eqref{base}
that $\mathbf{Q}_2$ is a coercive quadratic functional, which
guarantees $\mathbf{Q}_1$ is also coercive. Thus $\mathbf{Q}$ is a
coercive quadratic functional, which yields that the Lagrange
functional \eqref{ma} subject to \eqref{ma1} and \eqref{ma2} has a
minimum value at the unique critical point $(i_3^*,i_5^*,i_6^*)\in
L^2(t_0,t_1;\mathbb{R}^3)$. Following the approach to indefinite
linear quadratic optimal control problems in \cite{WL}, we can prove
any optimal control are continuous, which yields that
$(i_3^*,i_5^*,i_6^*)\in C([t_0,t_1];\mathbb{R}^3)$.

{\bf (II)} Let
\begin{equation}\label{construction}
\begin{array}{rcl}
\begin{cases}
\widehat{i_3}=\sum_{n=1}^{+\infty}\widehat{b_{3,n}}\widetilde{e_n}:=h\sum_{n=1}^{+\infty}\frac{4\pi^2C_1n}{4\pi^2C_1K(C_1)n^2-(t_1-t_0)^2}\widetilde{e_n},\\
\widehat{i_5}=\sum_{n=1}^{+\infty}\widehat{b_{5,n}}\widetilde{e_n}:=h\sum_{n=1}^{+\infty}\frac{2\pi^2C_2n}{2\pi^2C_2\widetilde{K}(C_2)n^2-(t_1-t_0)^2}\widetilde{e_n},\\
\widehat{i_6}=\sum_{n=1}^{+\infty}\widehat{b_{6,n}}\widetilde{e_n}:=h\sum_{n=1}^{+\infty}\frac{2\pi^2C_2n}{2\pi^2C_2\widetilde{K}(C_2)n^2-(t_1-t_0)^2}\widetilde{e_n},
\end{cases}
\end{array}
\qquad h\in\mathbb{R}.
\end{equation}

Lemma 6.2 and 6.3 yields that
\begin{eqnarray*}
\quad\begin{array}{rl}
&\frac{4}{t_1-t_0}\mathbf{Q}|_{(\widehat{i_3},\widehat{i_5},\widehat{i_6})}\\
=&h^2\{L_4\sum_{n=1}^{+\infty}(\widehat{b_{3,n}}-\widehat{b_{5,n}}-\widehat{b_{6,n}})^2\\
&+[L_3-K(C_1)]\sum_{n=1}^{+\infty}\widehat{b_{3,n}}^2+[L_5-\widetilde{K}(C_2)]\sum_{n=1}^{+\infty}\widehat{b_{5,n}}^2+[L_6-\widetilde{K}(C_2)]\sum_{n=1}^{+\infty}\widehat{b_{6,n}}^2\}\\
\le&h^2\{L_4\sum_{n=1}^{+\infty}(\widehat{b_{3,n}}-\widehat{b_{5,n}}-\widehat{b_{6,n}})^2\\
&+[L_3-K_2]\sum_{n=1}^{+\infty}\widehat{b_{3,n}}^2+[L_5-2K_2]\sum_{n=1}^{+\infty}\widehat{b_{5,n}}^2+[L_6-2K_2]\sum_{n=1}^{+\infty}\widehat{b_{6,n}}^2\}\\
=&h^2\sum_{n=1}^{+\infty}(\widehat{b_{3,n}},\widehat{b_{5,n}},\widehat{b_{6,n}})S_2(\widehat{b_{3,n}},\widehat{b_{5,n}},\widehat{b_{6,n}})^T<0,\end{array}
\end{eqnarray*}
which implies
$\lim_{h\rightarrow\infty}J(\widehat{i_3},\widehat{i_5},\widehat{i_6})=-\infty$.
So it follows from the density of $C([t_0,t_1];\mathbb{R}^3)$ in
$L^2(t_0,t_1;\mathbb{R}^3)$ that the Lagrange functional \eqref{ma}
subject to \eqref{ma1} and \eqref{ma2} has no minimum value.

{\bf (III)} If $C_1=C_2$, then $\widetilde{K}(C_2)=2K(C_1)$,
$K_1=K_2$ and $S_1=S_2$. If $S_2$ has at least one negative
characteristic root, then there exist at least one unit vector
$(x,y,z)$ such that
$$(x,y,z)S_2(x,y,z)^T<0.$$

Let
\begin{equation}\label{construction}
\begin{array}{rcl}
\begin{cases}
\widehat{i_3}=\sum_{n=1}^{+\infty}\widehat{b_{3,n}}\widetilde{e_n}:=hx\sum_{n=1}^{+\infty}\frac{4\pi^2C_1n}{4\pi^2C_1K(C_1)n^2-(t_1-t_0)^2}\widetilde{e_n},\\
\widehat{i_5}=\sum_{n=1}^{+\infty}\widehat{b_{5,n}}\widetilde{e_n}:=hy\sum_{n=1}^{+\infty}\frac{4\pi^2C_1n}{4\pi^2C_1K(C_1)n^2-(t_1-t_0)^2}\widetilde{e_n},\\
\widehat{i_6}=\sum_{n=1}^{+\infty}\widehat{b_{6,n}}\widetilde{e_n}:=hz\sum_{n=1}^{+\infty}\frac{4\pi^2C_1n}{4\pi^2C_1K(C_1)n^2-(t_1-t_0)^2}\widetilde{e_n},
\end{cases}
\end{array}
\qquad h\in\mathbb{R}.
\end{equation}
Analogous to the proof of {\bf (II)}, Lemma 6.2 and 6.3 yields that
\begin{eqnarray*}
\frac{4}{t_1-t_0}\mathbf{Q}|_{(\widehat{i_3},\widehat{i_5},\widehat{i_6})}
=h^2(x,y,z)S_2(x,y,z)^T\sum_{n=1}^{+\infty}\frac{16\pi^4C_1^2n^2}{[4\pi^2C_1K(C_1)n^2-(t_1-t_0)^2]^2}<0,
\end{eqnarray*}
which implies
$\lim_{h\rightarrow\infty}J(\widehat{i_3},\widehat{i_5},\widehat{i_6})=-\infty$.
So the Lagrange functional \eqref{ma} subject to \eqref{ma1} and
\eqref{ma2} has no minimum value.  \hfill$\Box$\vspace{3mm}

{\bf Proof of Proposition 2.2}\quad By Fubini Theorem, it follows
from \eqref{ma} and \eqref{ma1} that
\begin{eqnarray}\label{var}
\begin{array}{ll}
&J(i_3+\varepsilon\delta i_3,i_5+\varepsilon\delta i_5,i_6+\varepsilon\delta i_6)-J(i_3,i_5,i_6)\\
=&\varepsilon\int_{t_0}^{t_1}\{[L_3i_3+L_4(i_3-i_5-i_6)]\delta
i_3-\frac{1}{C_1}(q_1(t_0)+\int_{t_0}^ti_3\,d\tau)\int_{t_0}^t\delta
i_3\,ds\}\,dt\\
&+\varepsilon\int_{t_0}^{t_1}[L_5i_5+L_4(i_5+i_6-i_3)]\delta
i_5\,dt\\&\qquad\qquad\qquad-\varepsilon\int_{t_0}^{t_1}\frac{1}{C_2}[q_2(t_0)+\int_{t_0}^t(i_5+i_6)\,d\tau]\int_{t_0}^t\delta
i_5\,ds\,dt\\
&+\varepsilon\int_{t_0}^{t_1}[L_6i_6+L_4(i_5+i_6-i_3)]\delta
i_6\,dt\\&\qquad\qquad\qquad-\varepsilon\int_{t_0}^{t_1}\frac{1}{C_2}[q_2(t_0)+\int_{t_0}^t(i_5+i_6)\,d\tau]\int_{t_0}^t\delta
i_6\,ds\,dt+o(\varepsilon)\\
=&\varepsilon\int_{t_0}^{t_1}[L_3i_3+L_4(i_3-i_5-i_6)-\frac{1}{C_1}\int_t^{t_1}(q_1(t_0)+\int_{t_0}^si_3\,d\tau)\,ds]\delta
i_3\,dt\\
&+\varepsilon\int_{t_0}^{t_1}\{L_5i_5+L_4(i_5+i_6-i_3)-\frac{1}{C_2}\int_t^{t_0}[q_2(t_0)+\int_{t_0}^s(i_5+i_6)\,d\tau]\,ds\}\delta
i_5\,dt\\
&+\varepsilon\int_{t_0}^{t_1}\{L_6i_6+L_4(i_5+i_6-i_3)-\frac{1}{C_2}\int_t^{t_0}[q_2(t_0)+\int_{t_0}^s(i_5+i_6)\,d\tau]\,ds\}\delta
i_6\,dt\\
&+o(\varepsilon).
\end{array}
\end{eqnarray}

The constraints \eqref{ma2} yields that
\begin{eqnarray}\label{in}
\delta i_3\in\mathbb{H}_1,\quad \delta i_5\in\mathbb{H}_1,\quad
\delta i_6\in\mathbb{H}_1,
\end{eqnarray}
where
$$\mathbb{H}_1:=\{\sum_{n=1}^{+\infty}(a_ne_n+b_n\widetilde{e_n})|\hskip
1mm \sum_{n=1}^{+\infty}(a_n^2+b_n^2)<+\infty\},$$ and
$\mathbb{H}_0:=\{a e_0|\hskip 1mm a\in\mathbb{R}\}$, then
$L^2(t_0,t_1;\mathbb{R})=\mathbb{H}_0\oplus\mathbb{H}_1$, i.e.,
$\mathbb{H}_0$ is the orthogonal complement space of $\mathbb{H}_1$
in $L^2(t_0,t_1;\mathbb{R})$.

Hence, we have from \eqref{var}, \eqref{in} and
$L^2(t_0,t_1;\mathbb{R})=\mathbb{H}_0\oplus\mathbb{H}_1$ that there
exist some $l_3,l_5,l_6\in\mathbb{R}$ such that $(i_3,i_5,i_6)$
satisfy \eqref{ma2} and
\begin{equation}\label{sn}
\begin{array}{rcl}
\begin{cases}
L_3i_3+L_4(i_3-i_5-i_6)-\frac{1}{C_1}\int_t^{t_1}[q_1(t_0)+\int_{t_0}^si_3\,d\tau]\,ds=l_3, \\
L_5i_5+L_4(i_5+i_6-i_3)-\frac{1}{C_2}\int_t^{t_0}[q_2(t_0)+\int_{t_0}^s(i_5+i_6)\,d\tau]\,ds=l_5,\\
L_6i_6+L_4(i_5+i_6-i_3)-\frac{1}{C_2}\int_t^{t_0}[q_2(t_0)+\int_{t_0}^s(i_5+i_6)\,d\tau]\,ds=l_6,
\end{cases}
\end{array}
\end{equation}
if and only if $(i_3,i_5,i_6)$ is a critical point for the Lagrange
functional \eqref{ma} subject to \eqref{ma1} and \eqref{ma2}.
Through setting
\begin{equation*}
x_1(t):=\int_{t_0}^ti_3\,d\tau,\quad
x_2(t):=\int_{t_0}^ti_5\,d\tau,\quad x_3(t):=\int_{t_0}^ti_6\,d\tau,
\end{equation*}
it follows from the equation \eqref{sn} and \eqref{ma2} that
\begin{equation}\label{zong}
\begin{array}{rcl}
\begin{cases}
(L_4+L_3)x_1''-L_4x_2''-L_4x_3''+\frac{1}{C_1}x_1+\frac{q_1(t_0)}{C_1}=0, \\
-L_4x_1''+(L_4+L_5)x_2''+L_4x_3''+\frac{1}{C_2}x_2+\frac{1}{C_2}x_3+\frac{q_2(t_0)}{C_2}=0,\\
-L_4x_1''+L_4x_2''+(L_4+L_6)x_3''+\frac{1}{C_2}x_2+\frac{1}{C_2}x_3+\frac{q_2(t_0)}{C_2}=0.
\end{cases}
\end{array}
\end{equation}
with the boundary condition
\begin{equation}\label{zong1}
\begin{array}{rrr}
\begin{cases}
x_1(t_0)=0, \  x_2(t_0)=0, \  x_3(t_0)=0,\\x_1(t_1)=\lambda_3, \
 x_2(t_1)=\lambda_5, \  x_3(t_1)=\lambda_6.
\end{cases}
\end{array}
\end{equation}

Through defining
$$
y=(x_1,x_2,x_3,x_1',x_2',x_3')^T,
$$
and due to the positive definiteness of $M$, the boundary problem
\eqref{zong}-\eqref{zong1} can be reformulated as follows:
\begin{eqnarray}\label{zong2}
y'=\left(\begin{array}{ccc}0&I_3\\-M^{-1}N&0
\end{array}\right)y+\left(\begin{array}{ccc}0\\M^{-1}\mathbf{a}
\end{array}\right),
\end{eqnarray}
with the boundary condition
\begin{equation}\label{zong3}
y(t_0)=\left(\begin{array}{ccc}0\\\mathbf{c}
\end{array}\right),\qquad y(t_1)=\left(\begin{array}{ccc}x(t_1)\\\mathbf{d}
\end{array}\right),
\end{equation}
where $I_3$ is the $3\times 3$ identity matrix,
$x(t_1)=(\lambda_3,\lambda_5,\lambda_6)^T$ and the matrice $M$, $N$
and $\mathbf{a}$ are defined by \eqref{matrix}, \eqref{matrix1} and
\eqref{matrix4} while $\mathbf{c},\mathbf{d}\in\mathbb{R}^3$ are to
be known.

By the variation-of-constants formula, the problem
\eqref{zong2}-\eqref{zong3} is equivalent to
\begin{equation}\label{zong4}
\left(\begin{array}{ccc}\Phi(t_1-t_0)\mathbf{c}\\\Psi(t_1-t_0)\mathbf{c}
\end{array}\right)+\int_{t_0}^{t_1}\left(\begin{array}{ccc}\Phi(t_1-t)M^{-1}\mathbf{a}\\\Psi(t_1-t)M^{-1}\mathbf{a}
\end{array}\right)\,dt=\left(\begin{array}{ccc}x(t_1)\\\mathbf{d}
\end{array}\right),
\end{equation}
with
\begin{equation*}
\begin{array}{rcl}
\begin{cases}
\Psi(t):=I_3+\sum_{k=1}^{+\infty}\frac{t^{2k}}{(2k)!}(-M^{-1}N)^k,\\
\Phi(t):=t[I_3+\sum_{k=1}^{+\infty}\frac{t^{2k}}{(2k+1)!}(-M^{-1}N)^k].
\end{cases}
\end{array}
\end{equation*}

In the second equation of \eqref{zong4}, $\mathbf{d}$ is uniquely
determined by $\mathbf{c}$. So we only need to consider the
solvability of $\mathbf{c}$ through the first equation of
\eqref{zong4}.

By the definition of the matrice $M^{-\frac{1}{2}}$,
$M^{-\frac{1}{2}}NM^{-\frac{1}{2}}$ and $P$ in \eqref{matrix2} and
\eqref{matrix3}, it follows from the definition of $\Phi(t)$ that
\begin{equation}
\begin{array}{ll}
\frac{1}{t}PM^{\frac{1}{2}}\Phi(t)M^{-\frac{1}{2}}P^T&=I_3+\sum_{k=1}^{+\infty}\frac{t^{2k}}{(2k+1)!}(-PM^{-\frac{1}{2}}NM^{-\frac{1}{2}}P^T)^k\\
&=\left(\begin{array}{ccc}
\frac{1}{\sqrt{h_1}t}\sin(\sqrt{h_1}t)&0&0\\
0&\frac{1}{\sqrt{h_2}t}\sin(\sqrt{h_2}t)&0\\
0&0&1
\end{array}\right).
\end{array}
\end{equation}
Thus we obtain (I) from the invertibility of $\Phi(t_1-t_0)$. The
proof of (II) can be obtained by direct calculations in this case
$\Phi(t_1-t_0)$ is a singular matrix. \hfill$\Box$

\end{document}